\documentclass[12pt]{amsart}
\usepackage{amsmath, amssymb}

\author{}
\title{}
\date{}
\newcommand{\bbC}{{\mathbb C}}

\newcommand{\bbQ}{{\mathbb Q}}
\newcommand{\bbR}{{\mathbb R}}
\newcommand{\bbZ}{{\mathbb Z}}

\newcommand{\al}{\alpha}
\newcommand{\be}{\beta}
\newcommand{\hgt}{{\rm h}}

\newtheorem{lemma}{Lemma}
\newtheorem{theorem}{Theorem}

\begin{document}
\baselineskip=17pt

\vspace{10mm} 

\begin{center}
{\bf Primitive divisors of Lucas and Lehmer sequences, III} \\ 
\vspace{7 mm}
{\sc By} PAUL M VOUTIER \\

\vspace{3.0mm}

{\it Department of Mathematics, City University, London, EC1V 0HB, UK}
\end{center}

\vspace{3.0mm}

\begin{center}
{\it Abstract}
\end{center}

In this paper we prove that if $n > 30 \, 030$ then the 
$n$-th element of any Lucas or Lehmer sequence has a 
primitive divisor. 

\vspace{3.0mm}

\begin{center}
\rule{5.0cm}{0.3mm}
\end{center}

\begin{center}
1. {\it Introduction}
\end{center}

\vspace{3.0mm}

This is the third in a series of articles by the author attacking 
a conjecture that for $n > 30$, the $n$-th element of any Lucas or 
Lehmer sequence always has a primitive divisor. 

We start by recalling the definitions of these terms. 

Let $\al$ and $\be$ be algebraic numbers such that
$\al+\be$ and $\al \be$ are relatively prime non-zero 
rational integers and $\al / \be$ is not a root of unity. 
The sequence ${ \left( u_{n} \right) }_{n=0}^{\infty}$ 
defined by $\displaystyle u_{n} 
= \frac{\al^{n} - \be^{n}}{\al - \be}$ 
for $n \geq 0$ is called a {\it Lucas sequence}. 

If, instead of supposing that $\al+\be \in \bbZ$, we only 
suppose that $(\al+\be)^{2}$ is a non-zero rational integer, 
still relatively prime to $\al \be$, then we define the 
{\it Lehmer sequence} ${ \left( u_{n} \right) }_{n=0}^{\infty}$ 
associated to $\al$ and $\be$ by 
\begin{displaymath}
u_{n} = \displaystyle \frac{\al^{n}-\be^{n}}
			   {\al^{\delta} - \be^{\delta}},   
\end{displaymath}
where $\delta=1$ for odd $n$ and $\delta=2$ for even $n$. 

We say that a prime number $p$ is a {\it primitive divisor}
of a Lucas number $u_{n}$ if $p$ divides $u_{n}$ but does not
divide $(\al - \be)^{2}u_{2} \cdots u_{n-1}$. Similarly, $p$ is 
a primitive divisor of a Lehmer number $u_{n}$ if $p$ divides 
$u_{n}$ but not ${\left( \al^{2} - \be^{2} \right) }^{2}
u_{3} \cdots u_{n-1}$.

In the first article \cite{Vout1} of this series, we enumerated all 
Lucas and Lehmer sequences whose $n$-th element has no primitive 
divisor for certain $n \leq 12$ and all $12 < n \leq 30$ as well 
as presenting some strong numerical evidence to support the above 
conjecture. Then, in \cite{Vout2}, we established that the conjecture 
was true if the absolute logarithmic height of $\be/\al$, 
$\hgt (\be/\al)$, is at most four. 

Here we wish to determine a number $C$ such that if $n > C$ then 
$u_{n}$ always has a primitive divisor. For $\al$ and $\be$ real, 
the work of Carmichael \cite{Carm}, Ward \cite{Ward} and 
Durst \cite{Durst} shows that one can take $C=12$. Since the 
twelfth Fibonacci number has only $2$ and $3$, the third and 
fourth Fibonacci numbers, as its prime divisors, this result 
is best possible. 

When $\al$ and $\be$ are complex, we have $|\be/\al|=1$ which 
makes it much more difficult to obtain the necessary lower bounds 
for $\left| (\be/\al)^{d}-1 \right|$. Using lower bounds for 
linear forms in logarithms, Schinzel \cite{Sch} was the first 
to prove the existence of an absolute bound on $n$ in this case. 

Stewart \cite{Stew} established an explicit version of 
Schinzel's result, proving that one can take $C=e^{452}2^{67}$ 
for Lucas sequences and $C=e^{452}4^{67}$ for Lehmer sequences. 
In Theorem~2 of \cite{Vout2}, we showed that one may take 
$C=2 \cdot 10^{10}$ in both cases. The method of obtaining 
this result was the same as Stewart's, however, we benefited 
from the dramatic improvements in the lower bounds for linear 
forms in logarithms since the time of Stewart's work. 

In this article, we show that the value of $C$ may be taken 
even smaller.

\begin{theorem}
\label{thm:thm1}
For all $n > 30,030$, the $n$-th element of any  
Lucas or Lehmer sequence has a primitive divisor. 
\end{theorem}

Our principal tool in proving this theorem is a lower bound 
for linear forms in two logarithms. For many purposes, including 
our own here, the following result is the best presently known; 
its proof uses a judicious combination of the work of Laurent, 
Mignotte and Nesterenko \cite{LMN} and Liouville's inequality. 

\begin{theorem}
\label{thm:mylb}
{\rm (i)} Let $\gamma$ be a non-zero algebraic number of degree 
$d$ over $\bbQ$ with $d \hgt (\gamma) > 8$. Letting $\log \gamma$ 
denote the principal value of the logarithm of $\gamma$ and 
letting $b_{1}$ and $b_{2}$ be rational integers, we put 
\begin{displaymath}
\Lambda = b_{1} i \pi - b_{2} \log \gamma. 
\end{displaymath}

If $\Lambda \neq 0$ then 
\begin{displaymath}
\log |\Lambda| > -c_{1}(B) d^{3} \hgt (\gamma) \log^{2} B,  
\end{displaymath}
where 
\begin{displaymath}
c_{1}(B) = \min \left( \frac{\log 2}{16 d \log^{2} B} 
		       + \frac{B}{2d^{2}\log^{2} B}, 
		       22.36 { \left( \frac{1}{2} + \frac{0.28}{\log B} 
				      + \frac{1.8}{d \log B} \right) }^{2} 
		\right).
\end{displaymath}
and $B = \max \left( 2, \left| b_{1} \right|, \left| b_{2} \right| \right)$. 

\noindent
{\rm (ii)} Using the same notation and hypotheses as in $(i)$, 
if $\Lambda \neq 0$ then 
\begin{displaymath}
\log |\Lambda| > -9.1d^{3} \hgt (\gamma) \log^{2} B.  
\end{displaymath}
\end{theorem}

Our proof of Theorem~\ref{thm:thm1} proceeds in two steps. 
We use part~(ii) of Theorem~\ref{thm:mylb} to establish an  
initial improvement to our bound for $n$, while Theorem~\ref{thm:mylb}~(ii) 
will be used in direct calculations to eliminate smaller values 
of $n$. The integer $30030$ is the largest one for which we could 
not establish the necessary bounds needed to assert that $u_{n}$ 
always has a primitive divisor.

Theorem~\ref{thm:thm1} is quite an improvement on previous results,  
but the hardest part remains. One would need a significant 
decrease in the size of the constants in Theorem~\ref{thm:mylb} 
to bridge the gap completely. Quite likely, as is often the case 
when using lower bounds for linear forms in logarithms, a new 
approach will be necessary to complete the proof of the conjecture. 

\vspace{3.0mm}

\begin{center}
2. {\it Linear Forms in Logarithms} 
\end{center}

\vspace{3 mm}

To prove Theorem~\ref{thm:mylb} we will use Th\'{e}or\`{e}me~1 
of Laurent, Mignotte and Nesterenko \cite{LMN}. 

Let us first remind the reader of the notation of \cite{LMN}.  

We let $\gamma$ be an algebraic number of absolute value one, 
which is not a root of unity, with $\log \gamma$ as the principal 
value of its logarithm. Let $b_{1}$ and $b_{2}$ be positive 
integers and $\Lambda$ be as in the statement of 
Theorem~\ref{thm:mylb} and put $D = [\bbQ(\gamma):\bbQ]/2$. 

For any positive integer $S_{1}$, the set of numbers of the 
form $\pm \gamma^{s}$ with $0 \leq s < S_{1}$ has cardinality 
$2S_{1}$, so applying Th\'{e}or\`{e}me~1 of \cite{LMN} with 
$\alpha_{1}=-1, \log \alpha_{1}=\pi i, \alpha_{2}=\gamma, 
a_{1}=\rho \pi, a=a_{2}, R_{1}=2$ and $R_{2}=R-1$, we obtain 
the following result. 

\begin{lemma}
\label{lem:lmn1}
Let $K$ and $L$ be integers satisfying $K \geq 3$ and $L \geq 2$ 
and let $R,S_{1}$ and $S_{2}$ be positive integers with $R \geq 2$ 
and $2S_{1} \geq L$. Let $\rho > 1$. We put $N=KL, S=S_{1}+S_{2}-1$, 
\begin{displaymath}
g = \frac{1}{4}-\frac{N}{12RS} 
\hspace{3.0mm} \mbox{ and } \hspace{3.0mm}
b = ((R-1)b_{2}+(S-1)b_{1}) 
    { \left( \prod_{k=1}^{K-1} k! \right) }^{-2/(K^{2}-K)}.
\end{displaymath}

Let $a$ be a positive real number satisfying  
\begin{displaymath}
a \geq \rho |\log \gamma| + 2D \hgt (\gamma). 
\end{displaymath}

Suppose that 
\begin{equation}
\label{eq:zeroest}
\mbox{\rm Card}\ \{ rb_{2}+sb_{1}: 0 \leq r < R-1, 0 \leq s < S_{2} \} 
> (K-1)L,  
\end{equation}
and that 
\begin{equation}
\label{eq:cond}
K(L-1) \log \rho + (K-1) \log 2 - (D+1) \log N - D(K-1) \log b 
- gL(\rho \pi R + aS) > 0. 
\end{equation}

Then 
\begin{displaymath}
|\Lambda'| > \rho^{-KL+(1/2)}
\end{displaymath}
where 
\begin{displaymath}
\Lambda' = \Lambda \cdot \max \left\{ \frac{LSe^{LS|\Lambda|/(2b_{2})}}{2b_{2}},  
	   \frac{LRe^{LR|\Lambda|/(2b_{1})}}{2b_{1}} \right\}. 
\end{displaymath}
\end{lemma}

As a stepping stone to establishing Theorem~\ref{thm:mylb},  
we shall use Lemma~\ref{lem:lmn1} to first prove the following 
lemma which is essentially Th\'{e}or\`{e}me~3 of \cite{LMN} 
with some of the parameters and values changed in a manner 
more suited to our work here. 

\begin{lemma}
\label{lem:lmn3}
Let $b_{1}, b_{2}, \gamma, \log \gamma, \Lambda$ and $D$ 
be as above. Suppose further that $D \hgt (\gamma) > 4$. 
Also put $a = 6 \pi + 2D \hgt (\gamma)$ and 
\begin{displaymath}
H = \max \left( 6, D\log \left( \frac{b_{1}}{a}+\frac{b_{2}}{6\pi} \right)
		   + 2.95D + 1.8 \right). 
\end{displaymath}         

If $\Lambda \neq 0$ then 
\begin{displaymath}
\log |\Lambda| > -6.66aH^{2}.  
\end{displaymath}
\end{lemma}

\begin{center}
2.1. {\it Proof of Lemma} 2 
\end{center}

\vspace{3.0mm}

To prove this lemma, we will proceed as in Section~7 of \cite{LMN} 
where the authors prove their Th\'{e}or\`{e}me~3. 

We let $\rho = 6$ and $a = \rho \pi + 2D \hgt (\gamma)$. 
As in \cite{LMN}, we put 
\begin{eqnarray}
\label{eq:7.1}
\lambda & = & \log \rho, \nonumber \\ 
b'    & = & \frac{b_{1}}{a}+\frac{b_{2}}{\rho \pi}, \nonumber \\
h     & = & H-\lambda, \nonumber \\
L     & = & 2 + \left[ \frac{2h}{\lambda} \right], \nonumber \\
K     & = & 1 + \left[ k \rho \pi aL \right], \mbox{ where $k=0.15756$}, \\ 
R     & = & 2 + \left[ \sqrt{(K-1)La/(\rho \pi)} \right], \nonumber \\ 
S_{1} & = & \left[ \frac{L+1}{2} \right] \mbox{ and } \nonumber \\
S_{2} & = & 1+ \left[ \sqrt{\rho \pi(K-1)L/a} \right].   \nonumber 
\end{eqnarray}

As $H \geq 6 > 3\lambda$, from our choice of $L$, we have 
\begin{equation}
\label{eq:7.2}
6 \leq L \leq \frac{2H}{\lambda}. 
\end{equation}

Therefore, the condition $L \geq 2$ from Lemma~\ref{lem:lmn1} 
holds. Notice as well that the hypothesis $2S_{1} \geq L$ is 
satisfied and since $L \geq 2$, $S_{1}$ must be a positive 
integer, as required. 

Since $a > 26.8$, our expression for $K$ yields 
\begin{equation}
\label{eq:7.4} 
K \geq 79.5L \geq 477.
\end{equation}
So the requirement that $K \geq 3$ is certainly met.

Thus $(K-1)La/(\rho \pi)$ and $\rho \pi (K-1)L/a$ are both 
positive and as a consequence $R \geq 2$ and $S_{2} \geq 1$. 
Now all the conditions at the beginning of the statement 
of Lemma~\ref{lem:lmn1} are satisfied. 

From the definition of $H$, we deduce that 
\begin{equation}
\label{eq:7.5} 
h \geq D \log b' + 2.95D.
\end{equation}

Still following \cite{LMN}, we prove the following lemmas. 

\begin{lemma}
\label{lem:10}
Recalling that 
\begin{displaymath}
b = ((R-1)b_{2}+(S-1)b_{1}) 
    { \left( \prod_{k=1}^{K-1} k! \right) }^{-2/(K^{2}-K)}, 
\end{displaymath}
with the above definitions, we have 
\begin{displaymath}
\log b \leq \log (b') + 2.5 - \frac{\log(2\pi K/\sqrt{e})}{K-1}  
       \leq \frac{h}{D} - 0.45 - \frac{\log(2\pi K/\sqrt{e})}{K-1}.  
\end{displaymath}
\end{lemma}

\begin{proof}
From (\ref{eq:7.2}) and (\ref{eq:7.4}), it follows that 
$R-2 = \left[ \sqrt{(K-1)La/(\rho \pi)} \right] \geq 63$ and 
so we have $R-1 < 1.02(R-2)$. 

In addition,  
\begin{displaymath}
K-1 = \left[ k \rho \pi a L \right] \geq k \rho \pi a L - 1
= \frac{k \rho \pi a L - 1}{k \rho \pi a L} k \rho \pi a L. 
\end{displaymath}

From (\ref{eq:7.2}), we now see that $K-1 \geq (476/477)k \rho \pi aL$. 
Using this inequality and proceeding in a similar way, we have 
$S_{2}-1 \geq (43/44) \sqrt{(476k/477)} \rho \pi L$. Hence 
\begin{displaymath}
\frac{S_{1}-1}{S_{2}-1} \leq \frac{L-1}{0.775 \rho \pi L} < 0.07.  
\end{displaymath}

Since $(R-2) \rho \pi, \left( S_{2}-1 \right) a 
\leq \sqrt{\rho \pi a (K-1)L}$, it follows that 
\begin{eqnarray*}
(R-1)b_{2}+(S-1)b_{1} & \leq & 1.07((R-2)b_{2}+(S_{2}-1)b_{1}) \\  
		      & \leq & 1.07b'\sqrt{\rho \pi a(K-1)L}. 
\end{eqnarray*}                      

The remainder of the proof follows the proof of Lemme~10 in \cite{LMN}.   
\end{proof}

\vspace{3.0mm}

Our next lemma is a slight refinement of Lemme~11 of \cite{LMN}. 

\begin{lemma}
\label{lem:11}
Letting $g$ be as in Lemma~$\ref{lem:lmn1}$, we have 
\begin{displaymath}
gL(\rho \pi R + aS) \leq \left( \frac{\sqrt{k}}{3} + \frac{1}{6\rho\pi} 
			 \right) a \rho \pi L^{2} + 0.64aL. 
\end{displaymath}                       
\end{lemma}

\begin{proof}
The proof is identical to the proof of Lemme~11 in \cite{LMN}, 
except that our upper bound for $\sigma=aL(1/6+2\rho \pi/(3a))$ 
differs due to our choice of $\rho$ and $a$. 
\end{proof}

\vspace{7.0mm}

\begin{center}
2.1.1. {\it Examination of} (\ref{eq:cond}) 
\end{center}

\vspace{3.0mm}

We let 
\begin{displaymath}
\Phi = K(L-1)\lambda + (K-1) \log 2 - (D+1)\log(KL) - D(K-1) \log b 
       - gL(\rho \pi R + aS). 
\end{displaymath}       

This is the left-hand side of (\ref{eq:cond}). From Lemmas~\ref{lem:10} 
and \ref{lem:11}, it follows that 
\begin{eqnarray*}
\Phi & \geq & K(L-1) \lambda + (K-1) \log 2 - (D+1) \log (KL) - (K-1)h 
	      + 0.45(K-1)D \\ 
     &      & + D \log \left( \frac{2\pi K}{\sqrt{e}} \right) 
	      - \left( \frac{\sqrt{k}}{3} + \frac{1}{6\rho \pi} \right) 
			 \rho \pi a L^{2} - 0.64aL \\
     & \geq & \frac{K(L-1) \lambda}{2} + (K-1) \log 2 - \log K  
	      - (D+1) \log L + h + 0.45KD + 0.88D \\ 
     &      & - \left( \frac{\sqrt{k}}{3} + \frac{1}{6\rho \pi} \right)
			 \rho \pi a L^{2} - 0.64aL,  
\end{eqnarray*}                         
the last inequality following from collecting like terms together 
and noting that $\lambda(L-1)-h \geq \lambda(L-1)/2$ because of our 
choice of $L$. We now write the last expression as the sum of 
two terms $\Phi_{1}$ and $\Phi_{2}$, where 
\begin{eqnarray*}
\Phi_{1} & = & \frac{KL\lambda}{2} 
	       - \left( \frac{\sqrt{k}}{3} + \frac{1}{6\rho \pi} \right)
		 \rho \pi a L^{2} \mbox{ and } \\
\Phi_{2} & = & \left( 0.45D + \log 2 - \frac{\lambda}{2} \right) K + h 
	       - 0.64aL - \log K - (D+1) \log L + 0.88D - \log 2. 
\end{eqnarray*}

We therefore have $\Phi \geq \Phi_{1}+\Phi_{2}$. 

Since $K \geq \rho \pi akL$, we obtain 
\begin{displaymath}
\frac{\Phi_{1}}{\rho \pi aL^{2}} 
\geq \frac{\lambda k}{2} - \frac{\sqrt{k}}{3} - \frac{1}{6\rho \pi}.  
\end{displaymath}

We want the quantity on the right to be positive and it is this 
condition which leads us to our choice of $k$. 

For $\Phi_{2}$, we note that $aL \leq K/(k \rho \pi) < 0.337K$, so 
\begin{eqnarray*}
\Phi_{2} & \geq & (0.45D-0.42)K + h - \log (KL) - D \log L + 0.88D - \log 2 \\ 
	 & \geq & D \Phi_{21}+\Phi_{22},  
\end{eqnarray*}         
where 
\begin{displaymath}
\Phi_{21} = 0.005K - \log L + 0.88 \hspace{3.0mm} \mbox{ and } \hspace{3.0mm}
\Phi_{22} = 0.025K - \log (KL) - \log 2. 
\end{displaymath}

Now $0.013K \geq L$ and $K \geq 477$, so both of these quantities 
are positive. Therefore, (\ref{eq:cond}) holds for the choice of 
parameters given in the beginning of Section~2.1.

\begin{center}
2.1.2. {\it Conclusion of the proof} 
\end{center}

\vspace{3.0mm}

We will first use Lemma~\ref{lem:lmn1} to show that 
\begin{equation}
\label{eq:7.7}
\log |\Lambda'| > -6.65aH^{2}. 
\end{equation}

We have already shown that $K \geq 3, L \geq 2, R \geq 2, 
S_{1} \geq 1, S_{2} \geq 1, 2S_{1} \geq L$ and that 
(\ref{eq:cond}) holds. Moreover, since $| \log \gamma | \leq \pi$, 
the required condition on $a$ also holds, so it only remains 
to consider (\ref{eq:zeroest}).  

Continuing to follow \cite{LMN}, we break the proof into two 
parts. First, we suppose that the numbers of the form 
$rb_{2}+sb_{1}$ with $0 \leq r < R-1$ and $0 \leq s < S_{2}$ 
are pairwise distinct. By the definition of $R$ and $S_{2}$, 
$RS_{2}$ is greater than $(K-1)L$, so (\ref{eq:zeroest}) holds. 

In this case, all the hypotheses of Lemma~\ref{lem:lmn1} are 
satisfied and so 
\begin{displaymath}
\log |\Lambda'| \geq -\lambda KL + \frac{\lambda}{2}.  
\end{displaymath}

From (\ref{eq:7.2}) and $a > 26.8$, we have 
\begin{eqnarray*}
\lambda KL & \leq & \lambda (1 + \rho \pi akL)L 
= \left( 1 + \frac{1}{\rho \pi akL} \right) \lambda \rho \pi akL^{2} \\
& \leq & \left( 1 + \frac{1}{\rho \pi akL} \right) 
	 \frac{4k\rho \pi}{\lambda} a H^{2} < 6.65aH^{2}.
\end{eqnarray*}         

We now consider the case when the numbers of the form 
$rb_{2}+sb_{1}$ with $0 \leq r < R-1$ and $0 \leq s < S_{2}$ 
are not pairwise distinct. In this case there exist two 
relatively prime integers $r$ and $s$ with 
\begin{displaymath}
0 < |r| < R-1, \hspace{5.0mm} 0 < |s| < S_{2} \hspace{3.0mm} 
\mbox{ and } \hspace{3.0mm} rb_{2}+sb_{1} = 0. 
\end{displaymath}

We may assume that $b_{1}$ and $b_{2}$ are also relatively prime, 
for otherwise we can remove the common factor and get a better lower 
bound for $|\Lambda|$. Thus $r=b_{1}$ and $s=b_{2}$, and since 
$L \geq 2$, we have $LR/(2b_{1}) > 1$ and $LS/(2b_{2}) > 1$. If 
$c > 0$, then $xe^{cx} > 1$ for all $x > 1$. So 
\begin{displaymath}
\max \left\{ \frac{LSe^{LS|\Lambda|/(2b_{2})}}{2b_{2}},  
	     \frac{LRe^{LR|\Lambda|/(2b_{1})}}{2b_{1}} \right\} 
> 1,              
\end{displaymath}
and $|\Lambda'| > |\Lambda|$. 

Since $rb_{2}+sb_{1}=0$, we can write 
$\Lambda = s \log \gamma - r \log (-1)$. 
Liouville's inequality along with the definitions 
of $a, K$ and $S_{2}$ now yield
\begin{eqnarray*}
\log |\Lambda'| 
&   >  & \log |\Lambda| 
	 \geq -D \log 2 - \left( S_{2}-1 \right) D \hgt (\gamma) 
	 \geq -D \log 2 - \left( S_{2}-1 \right) a/2 \\
& \geq & -D \log 2 - \sqrt{k} \rho \pi aL/2. 
\end{eqnarray*}

From the inequalities $L \leq 2H/\lambda$ and $H \geq 6$, 
we can verify that 
\begin{displaymath}
\rho \pi a\sqrt{k} L < 8.4aH  < 1.4aH^{2}.  
\end{displaymath}

Our definition of $H$ and the assumption that $b_{1}$ and 
$b_{2}$ are positive show that $H > -D \log (6\pi)+2.95D$ 
so  $D \log 2 < 52H$. Since $a > 26.8$ and $H \geq 6$, 
$D\log 2 < 0.33aH^{2}$. Therefore $D\log 2 + \sqrt{k} \rho \pi aL/2 
< 1.03aH^{2}$. This shows that (\ref{eq:7.7}) also holds 
if (\ref{eq:zeroest}) is not true. Hence (\ref{eq:7.7}) 
holds in either case. 

We now want to deduce Lemma~\ref{lem:lmn3} from (\ref{eq:7.7}). 
We establish this by contradiction. Suppose that 
\begin{equation}
\label{eq:7.8}
\log |\Lambda| \leq -6.66aH^{2}. 
\end{equation}

From the lower bound in (\ref{eq:7.2}), it is easy to show 
that $2 < 0.04 \sqrt{k} aL, 1/2 < 0.02 \sqrt{k} \rho \pi L$ 
and $L/2 < 0.07 \sqrt{k} \rho \pi L$. Combining these 
inequalities with our choice of parameters in (\ref{eq:7.1}) 
and with the upper bound in (\ref{eq:7.2}), we have 
\begin{displaymath}
\max \left( \frac{LS}{2b_{2}}, \frac{LR}{2b_{1}} \right) 
\leq \max \left( \frac{LS}{2}, \frac{LR}{2} \right) 
\leq \max \left( \frac{1.1\sqrt{k} \rho \pi L^{2}}{2}, 
		 \frac{1.1 \sqrt{k}aL^{2}}{2} \right) 
< 0.5aH^{2}. 
\end{displaymath}

Since $e^{-x} < 1/x$, from (\ref{eq:7.8}) we have 
$|\Lambda| < 1/(6.66aH^{2})$, so  
\begin{displaymath}
\max \left( \frac{LSe^{LS|\Lambda|/(2b_{2})}}{2b_{2}}, 
	    \frac{LRe^{LR|\Lambda|/(2b_{1})}}{2b_{1}} \right) 
\leq 0.5e^{1/13.3}aH^{2} < 0.55aH^{2}. 
\end{displaymath}

Therefore, 
\begin{displaymath}
\log |\Lambda'| < -\left( 6.66 - \frac{\log(0.55aH^{2})}{aH^{2}} \right)  
		  aH^{2} 
		< -6.653aH^{2}, 
\end{displaymath}
since $a > 26.8$ and $H \geq 6$. But this contradicts 
(\ref{eq:7.7}). Therefore, (\ref{eq:7.8}) cannot be correct 
and so our lemma follows. 

\vspace{7.0mm}

\begin{center}
2.2. {\it Proof of Theorem}~\ref{thm:mylb}
\end{center}

\vspace{3 mm}

As in the proof of Lemma~2 of \cite{Vout2}, we can assume that 
$|\gamma|=1$ and that both $b_{1}$ and $b_{2}$ are non-zero.
From the first of these assumptions and since $\gamma$ is not 
a root of unity, we have $d \geq 2$. Moreover, since $\gamma 
\not\in \bbR$, $D=d/2$.

Finally, we can suppose that both $b_{1}$ and $b_{2}$ 
are positive. Otherwise, we can write 
$\Lambda = \pm \left( b_{1}' \log \gamma + b_{2}' i \pi \right)$, 
where $b_{1}'$ and $b_{2}'$ are both positive integers. However, 
in this case, we can replace $\gamma$ by $\gamma^{-1}$, and since 
we are taking the principal value of the logarithm, we have 
$\Lambda = \pm \left( -b_{1}' \log \gamma + b_{2}' i \pi \right)$,
so we can put $\Lambda$ in the form used in Lemma~\ref{lem:lmn3}. 

From Liouville's inequality see Exercise 3.4 of \cite{Wald}), we obtain
\begin{displaymath}
\left| \Lambda \right| \geq 2^{-D} \exp \left( -DB \hgt (\gamma) \right).
\end{displaymath}

Since $d \hgt (\gamma) > 8$, we see that 
\begin{displaymath}
\log |\Lambda| > 
- \left( \frac{\log 2}{16d \log^{2} B} + \frac{B}{2d^{2} \log^{2} B} \right) 
  d^{3} \hgt (\gamma) \log^{2} B.  
\end{displaymath}

To obtain the second term in the definition of $c_{1}(B)$, we 
turn to Lemma~\ref{lem:lmn3}. We first want to obtain an upper 
bound for $a$ in terms of $d$ and $\hgt (\gamma)$. We have  
\begin{displaymath}
a \leq d \hgt (\gamma)(6\pi/(d \hgt (\gamma))+1) 
< 3.357d \hgt (\gamma), 
\end{displaymath}
by our hypothesis that $d \hgt (\gamma) > 8$. 

Since $a > 26.8$, $1/a+1/(6 \pi) < 0.0904$, 
$H \leq \max(6,(d/2) \log B + 0.28d+1.8)$. Thus 
\begin{displaymath}
\log |\Lambda| 
> - 22.36 \max { \left( \frac{6}{d \log B}, 
			\frac{1}{2} + \frac{0.28}{\log B} 
			+ \frac{1.8}{d \log B} \right) }^{2} 
    d^{3} \hgt (\gamma) \log^{2} B.  
\end{displaymath}

Notice that if $6/(d \log B) \geq 1/2 + 0.28/\log B + 1.8/(d \log B)$, 
then $B \leq 38$ (since $d \geq 2$), in which case the Liouville 
estimate is stronger. Therefore, we can replace the $\max$ by just 
the second term in it.  

This establishes part (i) of Theorem~2. 

Part (ii) follows by noting that $c_{1}(B)$ is always less than 
9.1 for $d=2$, the maximum occurring near $B=5358$, and that, 
for fixed $B$, $c_{1}(B)$ decreases as $d$ increases.  

\vspace{7.0mm}

\begin{center}
3. {\it Preliminaries to the Proof of Theorem}~\ref{thm:thm1}
\end{center}

\vspace{3 mm}

Let $\al$ and $\be$ be algebraic numbers with $|\al| \geq |\be|$ 
and $\be/\al \not\in \bbR$ of degree $d$ satisfying $d \hgt (\be/\al) > 8$. 
We can prove the following lemma. 

\begin{lemma}
\label{lem:bounds}
{\rm (i)} If $n=1$ or $2$, then 
\begin{displaymath}
\log 2 + n\log |\al| \geq \log \left| \al^{n} - \be^{n} \right| 
\geq n \log |\al| - \left( \frac{d}{2} - 1 \right) \log 2 
     - \frac{d \hgt (\be/\al)}{2}.  
\end{displaymath}

\noindent
{\rm (ii)} For $n \geq 3$, we have  
\begin{displaymath}
\log 2 + n \log |\al| \geq \log \left| \al^{n} - \be^{n} \right| 
\geq n \log |\al| - c_{3} \left( n_{*} \right) d^{3} \hgt (\be/\al) 
		    \log^{2} n_{*}, 
\end{displaymath}
where $n_{*}=n/\gcd(n,2)$ and $c_{3}(n_{*})=c_{1}(n_{*})+0.02$.  
\end{lemma}

\begin{proof}
(i) The upper bound follows directly from the triangle 
inequality and our assumption that $|\al| \geq |\be|$. 

We can write $\al - \be = \al (1 - \be/\al)$. 
By Liouville's inequality along with the fact that $\be/\al$ 
is not a real number (see Exercise 3.4 of \cite{Wald}), 
\begin{displaymath}
\log \left| \be/\al \pm 1 \right| 
\geq - \left( \frac{d}{2} - 1 \right) \log 2 - \frac{d \hgt (\be/\al)}{2}.  
\end{displaymath}

The result for $n=1$ now follows. 

If $n=2$ then we can write $(\be/\al)^{2}-1$ as 
$(\be/\al-1)(\be/\al+1)$. Now one of these two 
terms is of absolute value greater than $\sqrt{2}$, 
so the result follows from the Liouville inequalities 
above.

(ii) The upper bound follows as before.  

For the lower bound, we consider the quantities 
$\log \left| (\be/\al)^{n_{*}} \pm 1 \right|$.  

Using the maximum modulus principle and Schwarz' lemma, 
if $z \in \bbC$ and $r \in \bbR$ with $0 < r < 1$ such 
that $|z| \leq r$, we have 
\begin{displaymath}
\left| \log (1 \pm z) \right| \leq \frac{|z \cdot \log (1-r)|}{r},  
\end{displaymath}

Applying this result for $1+z$ with $r=1/3$ and 
$z=e^{n_{*} \log (\be/\al)}-1$ where we take the principal 
value of the logarithm, we see that either 
\begin{displaymath}
\left| (\be/\al)^{n_{*}} - 1 \right| > 1/3 \mbox{ or } 
\left| \Lambda \right| = \left| n_{*} \log (\be/\al) - k \pi i \right| 
< 1.3 \left| (\be/\al)^{n_{*}} - 1 \right| < 0.5,  
\end{displaymath}
for some integer $k$.

In the first case, the lemma holds so we need only consider 
the second case. Since we took the principal value of 
the logarithm of $\be/\al$, in this case, we have 
$-\pi < \mbox{Im}(\log (\be/\al)) \leq \pi$ and so 
$|k| < n_{*}+0.5/\pi$ or $|k| \leq n_{*}$. 

As $\be/\al$ is, by assumption, not a root of unity, 
$\Lambda \neq 0$. Thus we may apply Theorem~\ref{thm:mylb} 
for $n_{*} \geq 2$, finding that 
\begin{displaymath}
\log \left| \Lambda \right| 
> - c_{1}(m) d^{3} \hgt (\be/\al) \log^{2} n_{*}. 
\end{displaymath}

For $n_{*} \geq 2$, we have 
$\log(1.3)/ \left( d^{3} \hgt (\be/\al) \log^{2} n_{*} \right) 
< 0.02$, upon recalling the hypotheses that $d \hgt (\be/\al) > 8$ 
and $d \geq 2$ (since $\be/\al \not\in \bbR$). Thus  
\begin{displaymath}
\log \left| \al^{n_{*}} - \be^{n_{*}} \right| 
\geq n_{*} \log |\al| - c_{3}(n_{*}) d^{3} \hgt (\be/\al) \log^{2} n_{*}.  
\end{displaymath}

By the same reasoning as above, we also have 
\begin{displaymath}
\log \left| \al^{n_{*}} + \be^{n_{*}} \right| 
\geq n_{*} \log |\al| - c_{3}(n_{*}) d^{3} \hgt (\be/\al) \log^{2} n_{*}. 
\end{displaymath}

Arguing as in part (i), we see that the lemma holds. 
\end{proof}

\vspace{3mm}

Finally we need bounds for the growth of arithmetic functions 
which will appear in the proof of Theorem~\ref{thm:thm1}. 

\begin{lemma}
\label{lem:arith1}
{\rm (i)} Let $\omega(n)$ denote the number of distinct 
prime factors of $n$. For $n \geq 3$, 
\begin{displaymath}
\omega(n) < \frac{1.3841 \log n}{\log \log n}.  
\end{displaymath}

\noindent
{\rm (ii)} For $n \geq 3$, 
\begin{displaymath}
\varphi(n) \geq \frac{n}{e^{\gamma}\log \log n + 2.50637/\log \log n}, 
\end{displaymath}
where $\gamma=0.57721\ldots$ is Euler's constant. 

\noindent
{\rm (iii)} Let $p_{1}=2,p_{2}=3,\ldots$ be the primes in increasing 
order. If $\omega(n) \geq 2$ then 
\begin{displaymath}
\varphi(n) \geq n \prod_{i=1}^{\omega(n)} \frac{p_{i}-1}{p_{i}}. 
\end{displaymath}
\end{lemma}

\begin{proof}
(i) This follows from Th\'{e}or\`{e}me~11 of \cite{Robin}. 

(ii) This is Theorem 15 of \cite{RS}. 

(iii) If $n=q_{1}^{a_{1}} \cdots q_{k}^{a_{k}}$ is the factorization 
of $n$ into primes then 
\begin{displaymath}
\varphi(n) = n \prod_{i=1}^{k} \frac{q_{i}-1}{q_{i}}. 
\end{displaymath}

Since $q_{i} \geq p_{i}$ for each $i$, the result follows. 
\end{proof}

\vspace{3 mm}

\begin{lemma}
\label{lem:arith2}
Suppose $n$ is a positive integer and let 
\begin{displaymath}
f(n) = \sum_{\stackrel{m|n}{\mu(m)=1}} \log^{2} (n/m). 
\end{displaymath}

\noindent
{\rm (i)} If $\omega(n) = 2$ then 
$f(n) < 2 \log^{2} n - 3.5 \log n + 3.5$. 

\noindent
{\rm (ii)} If $\omega(n) = 3$ then 
$f(n) < 4 \log^{2} n - 13.4 \log n + 15.5$. 

\noindent
{\rm (iii)} If $\omega(n) = 4$ then 
$f(n) < 8 \log^{2} n - 41 \log n + 71$. 

\noindent
{\rm (iv)} If $\omega(n) = 5$ then 
$f(n) < 16 \log^{2} n - 117 \log n + 264$. 

\noindent
{\rm (v)} If $\omega(n) = 6$ then 
$f(n) < 32 \log^{2} n - 309 \log n + 883$. 

\noindent
{\rm (vi)} If $\omega(n) = 7$ then 
$f(n) < 64 \log^{2} n - 775 \log n + 2718$. 

\noindent
{\rm (vii)} If $\omega(n) = 8$ then 
$f(n) < 128 \log^{2} n - 1886 \log n + 7913$. 

\noindent
{\rm (viii)} For any positive integer $n$, 
$f(n) < 2^{\omega(n)-1} \log^{2} n$.  
\end{lemma}

\begin{proof}
(i) It is easy to see that $f(n) \leq \log^2 (n) + \log^{2} (n/6)$, 
since 6 is the first integer with two distinct prime factors.  
Expanding this sum and rounding-off the numbers, we obtain the 
given inequality. 

(ii)-(vii) We proceed similarly in these cases. Suppose that 
$\omega(n)=k$. We take the sum of $\log^{2} n$ plus the sum 
of $\log^{2}(n/m)$ over the first ${k \choose 2}$ positive 
square-free integers $m$ with $\omega(m)=2$, plus the sum of 
$\log^{2}(n/m)$ over the first ${k \choose 4}$ positive 
square-free integers $m$ with $\omega(m)=4$, \ldots

(viii) This is just the trivial upper bound. There are 
$2^{\omega(n)-1}$ divisors $m$ of $n$ with $\mu(m)=1$ 
and for each of these $n/m \leq n$.  
\end{proof}

\vspace{7.0mm}

\begin{center}
4. {\it Proof of Theorem}~\ref{thm:thm1}
\end{center}

\vspace{3 mm}

We mentioned in the Introduction that if $\al$ and $\be$ 
are real then it is known that for $n > 12$, $u_{n}$ has  
a primitive divisor. Thus we may suppose that $\al$ and 
$\be$ are complex conjugates. 

Also notice that $\be/\al$ is a root of 
$\al \be X^{2} - \left( \al^{2}+\be^{2} \right) X +\al\be$.   
By our assumptions on $\al$ and $\be$, $\al^{2}+\be^{2} 
= (\al+\be)^{2}-2\al \be \in \bbZ$, so this polynomial has 
integer coefficients. Therefore, we can let $d=2$ in 
Lemma~\ref{lem:bounds} and we have $\hgt (\be/\al)=(\log|\al \be|)/2 
= \log | \al |$, since $\al$ and $\be$ are complex conjugates. 

It is well-known that 
\begin{displaymath}
\log \left| \Phi_{n}(\al,\be) \right|
= \sum_{m|n} \mu(m) \log \left| \al^{n/m} - \be^{n/m} \right| .   
\end{displaymath}

From the upper bounds in Lemma~\ref{lem:bounds} and the fact that 
\begin{displaymath}
\sum_{m|n} \mu(m) (n/m) = \varphi(n), 
\end{displaymath}
we see that 
\begin{equation}
\label{eq:philb}
\log \left| \Phi_{n}(\al,\be) \right|
\geq \varphi(n) \log |\al| - 2^{\omega(n)-1} \log 2 
     + \sum_{\stackrel{m|n}{\mu(m)=1}} 
       \log \left| { \left( \frac{\be}{\al} \right) }^{n/m} - 1 \right|. 
\end{equation}

Using the lower bounds in Lemma~\ref{lem:bounds} and our 
observation concerning the heights of $\al$ and $\be/\al$, 
we find that $\log \left| \Phi_{n} (\al,\be) \right|$ is  
greater than 
\begin{displaymath}
\varphi(n) \hgt (\be/\al) - 2^{\omega(n)-1} \log 2 - \hgt (\be/\al) 
- 8 \hgt (\be/\al) \sum_{\stackrel{m|n}{\mu(m)=1}} 
  c_{3} \left( { \left( \frac{n}{m} \right) }_{*} \right) 
  \log^{2} { \left( \frac{n}{m} \right) }_{*},  
\end{displaymath}
since $d=2$, the third term arising from the possibility 
that $\mu(n)$ or $\mu(n/2)$ is 1. 

Stewart \cite[Section~5]{Stew} has shown that if 
$\left| \Phi_{n}(\al,\be) \right| > n$ then $u_{n}$ 
has a primitive divisor. Combined with the fact,  
easily provable from Lemma~\ref{lem:arith1}(i) and 
some calculation, that 
$(\log n)/4 + 2^{\omega(n)-3} \log 2 + 1 < n^{1/4}$ 
for $n > 42$, this implies that we need only show that 
\begin{equation}
\label{eq:secondineq}
\varphi(n) - 8 \sum_{\stackrel{m|n}{\mu(m)=1}} 
	       c_{3} \left( { \left( \frac{n}{m} \right) }_{*} \right) 
	       \log^{2} { \left( \frac{n}{m} \right) }_{*} - n^{1/4} > 0, 
\end{equation}
since $\hgt (\be/\al) > 4$. 

Using Theorem~\ref{thm:mylb}(ii), we know that $c_{3}((n/m)_{*}) \leq 9.12$.  
Since $(n/m)_{*} \leq n/m$, from the bounds in Lemma~\ref{lem:arith1}(i) 
and (ii) along with Lemma~\ref{lem:arith2}(viii), we want to satisfy 
the inequality 
\begin{displaymath}
\frac{n}{e^{\gamma}\log \log n + 2.50637/\log \log n} 
- 36.5 \cdot 2^{1.3841 \log (n)/(\log \log n)} \log^{2} n 
- n^{1/4} > 0. 
\end{displaymath}

Using Maple~V, one can show that this inequality holds for 
$n > 18 \cdot 10^{6}$. Since $\omega(n) \geq 9$ implies that 
$n \geq 223,092,870$, we need only consider $\omega(n) \leq 8$.

If $n$ is prime then (\ref{eq:secondineq}) becomes 
$n - 1 - 9.12 \cdot 8 \cdot \log^{2}(n) - n^{1/4} > 0$.  
Again using Maple~V, we see that this inequality holds 
for $n > 5400$. 

Proceeding in a similar way, using Lemma~\ref{lem:arith1}(iii) 
and Lemma~\ref{lem:arith2}(i)--(vii), we find that $u_{n}$ 
has a primitive divisor for $\omega(n)=k$ if $n > n_{k}$ where 
the $n_{k}$'s are given in Table~1. 

\begin{table}[htb]
\centering
\begin{tabular}{||r|r|r|r|r|r|r|r|r||}                           \hline
$k$        &         1 &                  2 &                3 &         4 & 
		     5 &                  6 &                7 &         8    \\ \hline
$n_{k}$    &     5,400 &             43,000 &          109,000 &   256,000 &  
	       527,000 & $1.1 \cdot 10^{6}$ & $2 \cdot 10^{6}$ & $4 \cdot 10^{6}$ \\ \hline
\end{tabular}
\caption{Values of $n_{k}$}
\end{table}

Notice that $n_{8} < 9,699,690$ which is the smallest positive 
integer with eight distinct prime factors. Therefore, we do not 
need to consider the case of $\omega(n)=8$. 

If $\omega(n)=7$ then we need only look at $510,510 \leq n 
\leq 2 \cdot 10^{6}$. Notice that $3 \cdot 5 \cdots 19$, the 
smallest odd integer with seven distinct prime divisors, is 
greater than $n_{7}$. This implies that the $n$ in the given 
range with $\omega(n)=7$ are all even. Similarly, 
$2 \cdot 5 \cdot 7 \cdots 19 > n_{7}$ so these $n$ must also 
be divisible by 3. Continuing in this manner, we find that $n$ 
must be of the form 
$n=p_{1}^{a_{1}} p_{2}^{a_{2}} p_{3} p_{4} p_{5} p_{6} p_{7}$ 
with $p_{1}=2, p_{2}=3, p_{3}=5$ or 7, $p_{4}=7$ or 11, 
$p_{5}=11,13$ or 17, $p_{6}=13,17,19$ or 23, $17 \leq p_{7} \leq 61$,  
and $1 \leq a_{1},a_{2} \leq 2$. It turns out there are 58 of such 
integers and for each of them we calculated the left-hand side of 
(\ref{eq:secondineq}) using Maple~V. This quantity was positive 
except for the four cases 510510, 570570, 690690 and 746130. 

To eliminate these four cases we notice that if $n$ is odd  
then $\Phi_{2n}(\al,\be)=\Phi_{n}(\al,-\be)$. Therefore, 
if $\left| \Phi_{n}(\al,-\be) \right| > 2n$ then 
$\left| \Phi_{2n}(\al,\be) \right| > 2n$ and $u_{2n}$ has 
a primitive divisor.  

Some computation, along with Lemma~\ref{lem:arith1}(i), shows  
that $(\log (2n))/4 + 2^{\omega(n)-3} \log 2$  $+ 1 < n^{1/4}$ 
for $n > 210$, so (\ref{eq:secondineq}) implies that 
$|\Phi_{n}(\al,-\be)| > 2n$ if $n > 210$. We calculated the 
left-hand side of (\ref{eq:secondineq}) for $n=255255, 285285,
345345$ and 373065. In each case, it is positive. Therefore, 
$u_{510510}, u_{570570}, u_{690690}$ and $u_{746130}$ always 
have primitive divisors. We can conclude that if $\omega(n) = 7$ 
then $u_{n}$ has a primitive divisor.

We now consider the case of $\omega(n)=6$. For such 
$n=p_{1}^{a_{1}}p_{2}^{a_{2}}p_{3}^{a_{3}}p_{4}^{a_{4}}
p_{5}^{a_{5}}p_{6}^{a_{6}}$ with $n \leq 1.1 \cdot 10^{6}$, 
we have $2 \leq p_{1} \leq 3, 3 \leq p_{2} \leq 7, 
5 \leq p_{3} \leq 13, 7 \leq p_{4} \leq 29, 11 \leq p_{5} \leq 71, 
13 \leq p_{6} \leq 467, 1 \leq a_{1} \leq 6, 1 \leq a_{2} \leq 4, 
1 \leq a_{3} \leq 3$ and $1 \leq a_{4},a_{5},a_{6} \leq 2$. 
There are 2672 such positive integers and, again, we calculated 
the left-hand side of (\ref{eq:secondineq}) for each of these 
integers, this time using UBASIC~8.74 for the sake of speed. The 
largest two integers with $\omega(n)=6$ for which the left-hand 
side of (\ref{eq:secondineq}) is not positive are 480480 and 482790. 
The left-hand side of (\ref{eq:secondineq}) is positive for 
$n=241395$ and hence, arguing as above, $u_{482790}$ always has 
a primitive divisor. However, 480480 is not congruent to $2 \bmod 4$ 
so we cannot apply the same argument here. But there is one more 
idea we can take advantage of. 

Among the terms $|(\be/\al)^{n/m}-1|$ with $\mu(m)=1$ there 
will be one that is the smallest: suppose the minimum occurs 
for $m=r$. We can use this fact to complete the proof of the 
theorem. Let us write $(\be/\al)^{n/m}$ as $e^{i\theta_{m}}$ 
where $-\pi < \theta_{m} \leq \pi$. 

We first show that $|\theta_{r}| < 1/n^{2}$ for $n > 30$. If 
$|\theta_{r}| \geq 1/n^{2}$ then 
\begin{eqnarray*}
\log \left| (\be/\al)^{n/m}-1 \right| 
& \geq & \log \left| (\be/\al)^{n/r}-1 \right| 
	 = \log \sqrt{2-2\cos \theta_{r}} \\
&   >  & \log \sqrt{2-2\cos \left( 1/n^2 \right)} 
> -2.001\log n,
\end{eqnarray*}
for every $m$ which is a divisor of $n$ with $\mu(m)=1$ 
once $n > 30$. Using (\ref{eq:philb}) and our assumption 
that $\hgt (\al) > 4$ (recall that $d=2$), if 
\begin{displaymath}
\varphi(n) - 2^{\omega(n)-3}(\log (2) + 2.001 \log n) 
> \frac{\log n}{4}  
\end{displaymath}
then $u_{n}$ always has a primitive divisor.

Applying both parts of Lemma~\ref{lem:arith1}, one can show 
that this inequality is true if $n > 142$ and a quick direct 
calculation shows that the inequality holds for $30 < n < 143$ 
too. So we may assume that $|\theta_{r}| < 1/n^{2}$. 

We now want to find an expression for $\theta_{m}$ in terms of 
$\theta_{r}, r$ and $m$. Since $(\be/\al)^{n}=e^{ir\theta_{r}}
=e^{im\theta_{m}}$, we have 
\begin{displaymath}
\theta_{m} = \frac{r\theta_{r}}{m} + \frac{2\pi k}{m}, 
\end{displaymath}
for some $0 \leq k < m$.

Since $\theta_{r}$ is the minimum angle in absolute value, 
if $m > r$ then we must have $k \neq 0$. We saw that 
$|\theta_{r}| < 1/n^{2}$, so for $n > 30$, if $m > r$ 
then $|\theta_{m}| > 6/m$. From this we can conclude 
that for such $m$, 
\begin{displaymath}
\left| { \left( \frac{\be}{\al} \right) }^{n/m} - 1 \right| 
= \left| e^{i \theta_{m}} - 1 \right| 
= \sqrt{2-2\cos \theta_{m}} > \frac{1}{m}. 
\end{displaymath}

When $m \leq r$, we simply use  
\begin{displaymath}
\left| { \left( \frac{\be}{\al} \right) }^{n/m} - 1 \right| 
\geq \left| { \left( \frac{\be}{\al} \right) }^{n/r} - 1 \right|. 
\end{displaymath}

Returning to (\ref{eq:philb}) and applying Lemma~\ref{lem:bounds} 
only in the case where $n$ is replaced by $n/r$, we have 
\begin{eqnarray*}
\log | \Phi_{n}(\al,\be) | 
     & > & \hgt (\be/\al) F(n) \hspace{3.0mm} \mbox{ where } \\
F(n) & = & \min_{r|n,\mu(r)=1} 
	   \left( \varphi(n) - 2^{\omega(n)-3} \log 2 
		  - \frac{1}{4} \sum_{\stackrel{m|n, m > r}{\mu(m)=1}} 
		    \log m \right. \\
     &   & \left. \hspace*{20.0mm} 
		  - \sum_{\stackrel{m|n, m \leq r}{\mu(m)=1}} 
		  \max \left( 1,
		  8 c_{3} \left( (n/r)_{*} \right)
		  \log^{2} { \left( \frac{n}{r} \right) }_{*} \right)
		  \right), 
\end{eqnarray*}
since $\hgt (\be/\al) > 4$. The maximum expression, along 
with the $1$ inside it, in the last term are there to cover 
the possibility that $r=n/2$ or $r=n$.

Using UBASIC~8.74, we computed $F(n)$ for each of the 
$n \leq 480,480$ with $\omega(n)=6$ and (\ref{eq:secondineq}) 
not positive. We found that $F(n) > (\log n)/4$ for all 
$n > 265,650$. There are 61 values of $n$ less than or equal to 
265,650 with $\omega(n)=6$ and $F(n) \leq (\log n)/4$. Fortunately, 
each of these is congruent to 2 mod 4 and so $\log |\Phi_{n}(\al,\be)| 
> F(n/2)$. With the exception of $n=30030$ where $F(15015) = -691.61$, 
$F(n/2) > (\log n)/4$ for each of these $n$. As a consequence, 
Theorem~\ref{thm:thm1} holds if $\omega(n)=6$ and $n> 30030$. 
This whole procedure was automated using UBASIC~8.74 and its 
execution took about 200 seconds on an IBM-PC compatible computer 
with a 486DX CPU running at 33 MHz.

We proceed in the same way in the cases of $\omega(n)=3,4$ and 
$5$. There are $18569$ positive integers less than $527,000$ with 
$\omega(n)=5$. It took about 860 seconds to check, for all 
these values of $n$, whether $u_{n}$ always has a primitive 
divisor using $F(n)$ and $F(n/2)$, where appropriate. For all 
$n > 28980$, we were able to prove this while $F(28980)<-75.16$. 

We checked all 46412 positive integers with $\omega(n)=4$ and 
$n < 256,000$ in about 785 seconds. For all $n>26880$, either
$F(n)$, or $F(n/2)$ where appropriate, was greater than $\log(n)/4$  
and $F(26880)=-173.84$. 

Checking all the 42318 positive integers with $\omega(n)=3$ and 
$n < 109,000$ took about 445 seconds. For all $n>23040$, either
$F(n)$ or $F(n/2)$ was greater than $\log(n)/4$ and $F(23040)=-6.5$. 

Finally, we checked the 15491 positive integers with two distinct 
prime divisors which are less than 43,000 in about 100 seconds.  
Here we only used the simpler inequality (\ref{eq:secondineq}) 
and the largest such integer for which this inequality did not 
hold was 24,576, where the left-hand side was -94.09. 

These computations along with the results mentioned above 
comprise the proof of Theorem~\ref{thm:thm1}. 

To conclude, let us note that to ensure accuracy, in many cases, 
results obtained using UBASIC~8.74 were checked in Release~3 of 
Maple~V. This includes the computation of (\ref{eq:secondineq}), 
$F(n), F(n/2)$ and the number of integers with a given number of 
distinct prime divisors in certain intervals. In fact, the first 
three quantities were tested by hand for a few values of $n$. No 
discrepancy beyond round-off error was ever found. 

\setlength{\baselineskip}{5 mm}

\end{document}